\documentclass[a4paper,12pt]{article}
\usepackage{amsmath,amssymb,amsthm,amsxtra,mathrsfs,verbatim}
\newcommand{\C}{{\mathbb C}}
\newcommand{\N}{{\mathbb N}}
\newcommand{\R}{{\mathbb R}}

\newcommand{\abs}[2][\empty]{\ifx#1\empty\left|#2\right|%
\else#1\vert #2 #1\vert\fi}

\newcommand{\caninf}{\rho}
\newcommand{\Cnt}[1][]{{\cal C}^{#1}}
\newcommand{\co}[1]{{#1}^{c}}
\newcommand{\conv}{\star}
\newcommand{\csub}{\subset\subset}
\newcommand{\defstyle}[1]{{\em #1}}
\newcommand{\eps}{\varepsilon}

\renewcommand{\implies}{\Rightarrow}
\newcommand{\inv}[1]{{#1}^{-1}}

\newcommand{\monad}{\mu}
\newcommand{\norm}[2][\empty]{\ifx#1\empty\left\Vert#2\right\Vert%
\else#1\Vert #2 #1\Vert\fi}
\newcommand{\pairing}[2]{\langle#1,#2\rangle}
\newcommand{\restr}[2]{{#1}_{|#2}}
\newcommand{\Schwartz}{\mathscr{S}}
\newcommand{\sharpnorm}[2][\empty]{\abs[#1]{#2}_{\mathrm e}}

\newcommand{\supp}{\mathop{\mathrm{supp}}}
\newcommand{\test}{\mathcal{D}}

\newcommand{\val}{\mathrm{v}}

\newcommand{\Gen}{{\mathcal G}}
\newcommand{\GenC}{\widetilde\C}

\newcommand{\GenR}{\widetilde\R}
\newcommand{\Gentdinfty}{\mbox{\small \raisebox{-.2ex}{$\dot{\widetilde\Gen}$}}{}^\infty}
\newcommand{\E}{\mathcal{E}}
\newcommand{\EMod}{\E_M}
\newcommand{\Mod}{\mathcal{M}}
\newcommand{\Null}{\mathcal{N}}
\newcommand{\ns}{\mathrm{ns}}

\newtheorem{thm}{Theorem}[section]
\newtheorem{prop}[thm]{Proposition}
\newtheorem{lemma}[thm]{Lemma}
\newtheorem{cor}[thm]{Corollary}
\newtheorem{df}[thm]{Definition}
\newtheorem{ex}[thm]{Example}
\newtheorem{conj}[thm]{Conjecture}

\begin{document}
\title{Pointwise characterizations in generalized function algebras}
\author{Hans Vernaeve\footnote{Supported by research grants M949 and Y237 of the Austrian Science Fund (FWF)}\\
Unit for Engineering Mathematics\\
University of Innsbruck\\
A-6020 Innsbruck, Austria\\
{\tt Hans.Vernaeve@uibk.ac.at}}
\date{}
\maketitle
{
\abstract{We define the algebra $\Gen(\GenR^d)$ of Colombeau generalized functions on $\GenR^d$ which naturally contains the generalized function algebras $\Gen_\Schwartz(\R^d)$ and $\Gen_\tau(\R^d)$.
The subalgebra $\Gen_\Schwartz^\infty(\R^d)$ of $\Gen(\GenR^d)$ is characterized by a pointwise property of the generalized functions and their Fourier transforms.
We also characterize the equality in the sense of generalized tempered distributions for certain elements of $\Gen_\Schwartz(\R^d)$ (namely those with so-called slow scale support) by means of a pointwise property of their Fourier transforms.
Further, we show that (contrary to what has been claimed in the literature) for an open set $\Omega\subseteq\R^d$, the algebra of pointwise regular generalized functions $\dot\Gen^\infty(\Omega)$ equals $\Gen^\infty(\Omega)$ and give several characterizations of pointwise $\Gen^\infty$-regular generalized functions in $\Gen(\Omega)$.}
}



\section{Introduction}
Algebras of generalized functions have been developed by many authors \cite{DHPV04, GKOS, NPS} mainly inspired by the work of J.-F.\ Colombeau \cite{Colombeau85}, and have proved valuable as a tool for treating partial differential equations with singular data or coefficients (see \cite{O92} and the references therein). Under the influence of microlocal analysis of partial differential operators \cite{DPS, GGO05, GH05, HdH01, HOP05} and group invariance of solutions to partial differential equations \cite{Kunzinger97,Kunzinger00,Kunzinger06b} in the context of these algebras, the investigation of local properties of generalized functions became increasingly important. It was soon realized that generalized functions in the sense of Colombeau can be viewed naturally as pointwise functions on sets of generalized points \cite{KO99}. During the past years, there is a growing insight that many aspects of Colombeau generalized functions can be naturally expressed by pointwise properties (e.g., differentiability \cite{AFJ05}, regularity \cite{OPS}).\\
The purpose of this paper is to further investigate pointwise properties in Colombeau algebras. In section \ref{sec_GenGenRd}, we define the algebra $\Gen(\GenR^d)$ of generalized functions on $\GenR^d$, the space of generalized points of $\R^d$. In a sense, this is the largest algebra of generalized functions of Colombeau type that is embedded in the algebra of all pointwise maps $\GenR^d\to\GenC$ through the usual pointwise action of generalized functions.
Other known generalized function algebras such as $\Gen_c(\R^d)$, $\Gen_\Schwartz(\R^d)$ and $\Gen_\tau(\R^d)$ (see the next section for their definitions) are naturally contained in it and their elements can be characterized in $\Gen(\GenR^d)$ by a pointwise property (Theorems \ref{thm_E_tau_charac}, \ref{thm_E_Schwartz_charac}). 
We introduce the subalgebra $\Gen_{ss}(\GenR^d)$ of \defstyle{slow scale supported} generalized functions in $\Gen(\GenR^d)$, which can be characterized both by growth properties (similar to the definition of $\Gen_\Schwartz^\infty(\R^d)$ and of $\Gen_\Schwartz(\R^d)$) and by pointwise properties of the associated pointwise maps (Theorem \ref{thm_slowscale_support}). It is a subalgebra of $\Gen_\Schwartz(\R^d)$, and it turns out that $\Gen_{\Schwartz}^\infty(\R^d)=\Gen_{ss}(\GenR^d)\cap\widehat\Gen_{ss}(\GenR^d)$ (here $\widehat{\ }$ denotes the Fourier transform), which yields a pointwise characterization of elements in $\Gen_\Schwartz^\infty(\R^d)$ (Theorem \ref{thm_slowscale_spectrum}).\\
In section \ref{sec_gtd}, we characterize when elements of $\Gen_{ss}(\GenR^d)$ act identically on test functions in Schwartz's space of rapidly decreasing functions $\Schwartz(\R^d)$ (i.e., when they are equal in the sense of generalized tempered distributions) by means of a pointwise property of their Fourier transforms (Theorem \ref{thm_gds_characterization}).\\
Further, in section \ref{sec_regularity}, we discuss several notions of pointwise $\Gen^\infty$-regularity that were introduced in \cite{DPS,OPS} and show that they can all be characterized by means of one such notion ($\Gentdinfty$-regularity at a compactly supported generalized point), but for different sets of compactly supported generalized points (Corollary \ref{cor_Geninfty_charac}). Contrary to what has been claimed in the literature, we prove that for an open set $\Omega\subseteq\R^d$, the algebra of pointwise regular generalized functions $\dot\Gen^\infty(\Omega)$ in fact equals $\Gen^\infty(\Omega)$ (Theorem \ref{thm_Gendotinfty_equals_Geninfty_globally}) by means of a new characterization of $\dot\Gen^\infty$-regularity at a point (Corollary \ref{cor_Landau}). We conclude with a characterization of $\widetilde\Gen^\infty$-regularity at a point that might provide a basis for microlocal analysis of generalized functions on the scale of the sharp neighbourhoods (Proposition \ref{prop_Taylor}).

\section{Preliminaries}
In this paper, $\Omega$ denotes an open subset of $\R^d$ ($d\in\N$).\\
We recall the definitions of some generalized function algebras of Colombeau type \cite{Garetto05,GKOS,O92}.
the algebra of generalized functions on $\Omega$ equals $\Gen(\Omega)=\EMod(\Omega)/\Null(\Omega)$, where
\begin{align*}
\EMod(\Omega) =\,&\big\{(u_\eps)_\eps \in\Cnt[\infty](\Omega)^{(0,1)}: (\forall
K\csub\Omega) (\forall\alpha\in\N^d) (\exists N\in\N)\\
&(\exists\eps_0\in(0,1)) (\forall\eps\le\eps_0) \big(
\sup_{x\in K}\abs{\partial^\alpha u_\eps(x)}\le \eps^{-N}\big)\big\}\\
\Null(\Omega) =\,&\big\{(u_\eps)_\eps \in\Cnt[\infty](\Omega)^{(0,1)}: (\forall
K\csub\Omega) (\forall\alpha\in\N^d) (\forall m\in\N)\\
&(\exists\eps_0\in(0,1)) (\forall\eps\le\eps_0) \big(
\sup_{x\in K}\abs{\partial^\alpha u_\eps(x)}\le \eps^m\big)\big\}.
\end{align*}
The algebra $\Gen^\infty(\Omega)$ of regular generalized functions on $\Omega$ consists of those $u\in\Gen(\Omega)$ admitting a representative $(u_\eps)_\eps$ satisfying
\[
(\forall K\csub \Omega) (\exists N\in\N) (\forall\alpha\in\N^d)
(\exists\eps_0\in(0,1)) (\forall\eps\le\eps_0) \big(\sup_{x\in K}\abs[\big]{\partial^\alpha u_\eps(x)}\le\eps^{-N}\big).
\]
For a compact set $K\csub\R^d$, let $\test(K)$ denote the space of all $\Cnt[\infty](\R^d)$-functions with support contained in $K$.
The algebra $\Gen_c(\R^d)$ of compactly supported generalized functions on $\R^d$ consists of those $u\in\Gen(\R^d)$ admitting a representative $(u_\eps)_\eps\in\test(K)^{(0,1)}$, for some $K\csub\R^d$.\\
The algebra of tempered generalized functions on $\R^d$ equals $\Gen_\tau(\R^d)={\mathcal E}_\tau(\R^d)/\Null_\tau(\R^d)$, where
\begin{align*}
{\mathcal E}_{\tau}(\R^d)=\,&\big\{(u_\eps)_\eps\in\Cnt[\infty](\R^d)^{(0,1)}: (\forall\alpha\in\N^d)(\exists N\in\N)\\
&(\exists\eps_0\in(0,1)) (\forall\eps\le\eps_0) \big(\sup_{x\in\R^d}\abs[\big]{\partial^\alpha u_\eps(x)}(1 + \abs x)^{-N}\le\eps^{-N}\big)\big\}\\
{\Null}_{\tau}(\R^d)=\,&\big\{(u_\eps)_\eps\in\Cnt[\infty](\R^d)^{(0,1)}: (\forall\alpha\in\N^d)(\exists N\in\N) (\forall m\in\N)\\
&(\exists\eps_0\in(0,1)) (\forall\eps\le\eps_0) \big(\sup_{x\in\R^d}\abs[\big]{\partial^\alpha u_\eps(x)}(1 + \abs x)^{-N}\le\eps^{m}\big)\big\}.
\end{align*}
The algebra of generalized functions based upon Schwartz's space $\Schwartz(\R^d)$ of rapidly decreasing functions on $\R^d$ equals $\Gen_{\Schwartz}(\R^d)=\E_\Schwartz(\R^d)/\Null_\Schwartz(\R^d)$, where
\begin{align*}
{\E}_{\Schwartz}(\R^d)=\,&\big\{(u_\eps)_\eps\in\Schwartz(\R^d)^{(0,1)}: (\forall\alpha,\beta\in\N^d)(\exists N\in\N)\\
&(\exists\eps_0\in(0,1)) (\forall\eps\le\eps_0) \big(\sup_{x\in\R^d}\abs[\big]{x^\alpha\partial^\beta u_\eps(x)}\le\eps^{-N}\big)\big\}\\
{\Null}_{\Schwartz}(\R^d)=\,&\big\{(u_\eps)_\eps\in\Schwartz(\R^d)^{(0,1)}: (\forall\alpha,\beta\in\N^d) (\forall m\in\N)\\
&(\exists\eps_0\in(0,1)) (\forall\eps\le\eps_0) \big(\sup_{x\in\R^d}\abs[\big]{x^\alpha\partial^\beta u_\eps(x)}\le\eps^{m}\big)\big\}.
\end{align*}
The algebra $\Gen_\Schwartz^\infty(\R^d)$ of $\Schwartz$-regular generalized functions on $\R^d$ consists of those $u\in\Gen_\Schwartz(\R^d)$ admitting a representative $(u_\eps)_\eps$ satisfying
\[
(\exists N\in\N) (\forall\alpha,\beta\in\N^d)
(\exists\eps_0\in(0,1)) (\forall\eps\le\eps_0) \big(\sup_{x\in\R^d}\abs[\big]{x^\alpha\partial^\beta u_\eps(x)}\le\eps^{-N}\big).
\]
In all these algebras, the generalized function with representative $(u_\eps)_\eps$ is denoted by $[(u_\eps)_\eps]$.\\
We also recall the definition of generalized points and pointwise value theorems in these generalized function algebras.
The set of generalized points of $\R^d$ equals $\GenR^d=\Mod_{\R^d}/\Null_{\R^d}$, where
\begin{align*}
{\Mod}_{\R^d}=\,&\{(x_\eps)_\eps\in(\R^d)^{(0,1)}: (\exists N\in\N)
(\exists\eps_0\in(0,1)) (\forall\eps\le\eps_0) (\abs{x_\eps}\le\eps^{-N})\}\\
{\Null}_{\R^d}=\,&\{(x_\eps)_\eps\in(\R^d)^{(0,1)}: (\forall m\in\N)
(\exists\eps_0\in(0,1)) (\forall\eps\le\eps_0) (\abs{x_\eps}\le\eps^{m})\}.
\end{align*}
Similarly, the ring $\GenC$ of generalized (complex) numbers is defined using nets in $\C^{(0,1)}$.
By $\widetilde\Omega$, we denote the set of those $\tilde x\in\GenR^d$ admitting a representative $(x_\eps)_\eps\in\Omega^{(0,1)}$. By $\widetilde\Omega_c$, we denote the set of those $\tilde x\in\GenR^d$ admitting a representative $(x_\eps)_\eps\in K^{(0,1)}$ for some $K\csub\Omega$. The value $u(\tilde x)$ of a generalized function $u=[(u_\eps)_\eps]$ at a generalized point $\tilde x=[(x_\eps)_\eps]$ is defined as $[(u_\eps(x_\eps))_\eps]\in\GenC$. This definition turns out to be independent of representatives for $u\in\Gen(\Omega)$ at $\tilde x\in\widetilde\Omega_c$, and for $u\in\Gen_\tau(\R^d)$, $\Gen_\Schwartz(\R^d)$ at $\tilde x\in\GenR^d$.\\
Moreover, for $u\in\Gen(\Omega)$, $u=0$ iff $u(\tilde x)=0$ for each $\tilde x\in\widetilde\Omega_c$ \cite[Thm.\ 1.2.46]{GKOS}. Similarly, for $u\in\Gen_\tau(\R^d)$, $u=0$ iff $u(\tilde x)=0$ for each $\tilde x\in\GenR^d$ \cite[Prop.\ 1.2.47]{GKOS}. Since $\Null_\Schwartz(\R^d)=\Null_\tau(\R^d)\cap\E_\Schwartz(\R^d)$ \cite[Thm.\ 3.8]{Garetto05b}, the same holds for $u\in\Gen_\Schwartz(\R^d)$.\\
We denote by $\caninf\in\GenR$ the element $[(\eps)_\eps]$. For $\tilde x = [(x_\eps)_\eps]\in\GenR^d$, $\val(\tilde x)=\sup\{a\in\R: \abs{x_\eps}\le\eps^a$, for small $\eps\}$ denotes the valuation on $\GenR^d$ and $\sharpnorm{\tilde x}=e^{-\val(x)}$ the sharp norm on $\GenR^d$. The sharp norm defines a topology on $\GenR^d$, which is called the sharp topology (cf.\ \cite{Garetto05,Scarpalezos93}). For $x_0\in\R^d$, we denote by $\mathcal V_{x_0}$ the set of all (nongeneralized) neighbourhoods in $\R^d$ of $x_0$; for $\tilde x_0\in\GenR^d$, we denote by $\widetilde{\mathcal V}_{\widetilde x_0}$, the set of all sharp neighbourhoods in $\GenR^d$ of $\tilde x_0$. We refer to \cite{GKOS} for further properties of Colombeau generalized functions.  

\section{The algebra of generalized functions on $\GenR^d$}\label{sec_GenGenRd}
\begin{df}
$\Gen(\GenR^d):=\EMod(\GenR^d)/\Null(\GenR^d)$, where
\begin{align*}
\EMod(\GenR^d)=\,&\big\{(u_\eps)_\eps\in\Cnt[\infty](\R^d)^{(0,1)}: (\forall\alpha\in\N^d) (\forall m\in\N) (\exists N\in\N)\\ &(\exists\eps_0\in(0,1)) (\forall\eps\le\eps_0)
\big(\sup_{\abs{x}\le \eps^{-m}}\abs{\partial^\alpha u_\eps(x)}\le\eps^{-N}\big)\big\},\\
\Null(\GenR^d)=\,&\{(u_\eps)_\eps\in\Cnt[\infty](\R^d)^{(0,1)}: (\forall\alpha\in\N^d) (\forall m\in\N) (\forall n\in\N)\\
&(\exists\eps_0\in(0,1))(\forall\eps\le\eps_0)
\big(\sup_{\abs{x}\le \eps^{-m}} \abs{\partial^\alpha u_\eps(x)} \le\eps^n\big)\big\}.
\end{align*}
\end{df}
As for $\Gen(\R^d)$, it follows that $\Gen(\GenR^d)$ is a differential algebra and that we can always find representatives in $\test(\R^d)^{(0,1)}$.

\begin{prop}
\begin{multline*}
\Null(\GenR^d)=\EMod(\GenR^d)\cap \{(u_\eps)_\eps\in\Cnt[\infty](\R^d)^{(0,1)}: (\forall m\in\N) (\forall n\in\N)\\
(\exists\eps_0\in(0,1))(\forall\eps\le\eps_0)
\big(\sup_{\abs{x}\le \eps^{-m}} \abs{u_\eps(x)} \le\eps^n\big)\big\}.
\end{multline*}
Hence also in $\Gen(\GenR^d)$, $u=0$ iff $u(\tilde x)=0$, $\forall\tilde x\in\GenR^d$.
\end{prop}
\begin{proof}
As for $\Null(\R^d)$ \cite[Thm.\ 1.2.3, Thm.\ 1.2.46]{GKOS}.
\end{proof}

\begin{df}
Since $\EMod(\GenR^d)\subseteq\EMod(\R^d)$ and $\Null(\GenR^d)\subseteq\Null(\R^d)$, the identity map on representatives gives rise to a well-defined map $\Gen(\GenR^d)\to\Gen(\R^d)$. If $u\in\Gen(\GenR^d)$, we call the corresponding element $\restr{u}{\R^d}\in\Gen(\R^d)$ the \defstyle{restriction} of $u$ to $\Gen(\R)$.
\end{df}

The pointwise value theorems provide a natural way to see $\Gen_\Schwartz(\R^d)$ and $\Gen_\tau(\R^d)$ as subalgebras of $\Gen(\GenR^d)$. 
\begin{prop}\leavevmode
\begin{enumerate}
\item $\Gen_\tau(\R^d) = \mathcal{E}_\tau(\R^d)/(\mathcal{E}_\tau(\R^d)\cap\Null(\GenR^d))\subset\Gen(\GenR^d)$
\item $\Gen_\Schwartz(\R^d) = \mathcal{E}_\Schwartz(\R^d)/(\mathcal{E}_\Schwartz(\R^d)\cap\Null(\GenR^d)) \subset \Gen(\GenR^d)$
\end{enumerate}
\end{prop}
\begin{proof}
(1) As $\mathcal{E}_\tau(\R^d)\subseteq \EMod(\GenR^d)$, the point value characterizations in $\Gen_\tau(\R^d)$ and $\Gen(\GenR^d)$ imply that $\mathcal{E}_\tau(\R^d)\cap\Null(\GenR^d)=\Null_\tau(\R^d)$. Hence the identity map on representatives defines a canonical embedding of $\Gen_\tau(\R^d)$ into $\Gen(\GenR^d)$. Let $u_\eps(x)=(1+\abs x^2)^{\frac{\ln(1+\abs x^2)}{\ln(\eps^{-1})}}$, for each $\eps\in(0,1)$. As $\ln(1+\abs x^2)\le (2m+1)\ln(\eps^{-1})$, for each $x$ with $\abs x\le\eps^{-m}$ ($m\in\N$) and small $\eps$, the net $(u_\eps)_\eps$ is easily seen to belong to $\EMod(\GenR^d)$; as $u_\eps(x)\ge \eps^{-m^2}$, for each $x$ with $\abs{x}=\eps^{-m}$ ($m\in\N$), the net $(u_\eps)_\eps\notin \mathcal{E}_\tau(\R^d) + \Null(\GenR^d)$; hence the inclusion is strict.\\
(2) Similar.
\end{proof}
As in $\Gen(\R^d)$ (but in contrast with $\Gen_\tau(\R^d)$), we have a pointwise invertibility criterium in $\Gen(\GenR^d)$.
\begin{prop}
Let $u\in\Gen(\GenR^d)$. Then the following are equivalent:
\begin{enumerate}
\item there exists $v\in\Gen(\GenR^d)$ such that $uv=1$
\item for each $\tilde x\in\GenR^d$, $u(\tilde x)$ is invertible in $\GenC$
\item for some (and hence any) representative $(u_\eps)_\eps$ of $u$,
\[(\forall m\in\N) (\exists \eps_m \in (0,1)) (\exists n\in\N) (\forall\eps\le\eps_m) (\inf_{\abs x\le\eps^{-m}} \abs{u_\eps(x)} \ge \eps^n).\]
\end{enumerate}
\end{prop}
\begin{proof}
$(1)\Rightarrow(2)$: for $\tilde x\in\GenR^d$, $u(\tilde x)v(\tilde x)=1$ in $\GenC$.\\
$(2)\Rightarrow(3)$: supposing that the conclusion is not true, we find $M\in\N$ and a decreasing sequence $(\eps_n)_{n\in\N}$ tending to $0$ and $x_{\eps_n}\in\R^d$ with $\abs{x_{\eps_n}}\le\eps^{-M}$ and $\abs{u_{\eps_n}(x_{\eps_n})}<\eps_n^n$, for each $n\in\N$. Let $x_\eps=0$, if $\eps\notin\{\eps_n: n\in\N\}$. Then $\tilde x:=[(x_\eps)_\eps]\in\GenR^d$, but $u(\tilde x)$ is not invertible in $\GenC$ by \cite[Thm.\ 1.2.38]{GKOS}.\\
$(3)\Rightarrow(1)$: let $\eps_m$ as in the statement of (3). We may suppose that $(\eps_m)_{m\in\N}$ is decreasing and tends to $0$. Let $v_\eps\in\Cnt[\infty](\R^d)$ with $v_\eps(x)=\inv{u_\eps(x)}$, for $\abs x\le\eps^{-m}$ and $\eps_{m+1}<\eps\le\eps_m$. Since each $\partial^\alpha v_\eps$ is a linear combination (with coefficients indep.\ of $\eps$) of $\prod_{\beta}\partial^\beta u_\eps(x)/u_\eps^{\abs\alpha + 1}(x)$ (finite products), 
$[(v_\eps)_\eps]\in\Gen(\GenR^d)$. As $\sup_{\abs x\le \eps^{-m}}\abs{u_\eps(x) v_\eps(x) - 1} = 0$, for $\eps\le\eps_m$, we have $uv=1$ in $\Gen(\GenR^d)$.
\end{proof}

\begin{lemma}\label{lemma_eps_0_indep_of_m}
Let $u\in\Gen(\GenR^d)$. Then there exists a representative $(u_\eps)_\eps$ of $u$ satisfying
\[
(\forall\alpha\in\N^d) (\exists\eps_0\in(0,1)) (\forall m\in\N) (\exists N\in\N)\\
(\forall\eps\le\eps_0)\big(\sup_{\abs{x}\le \eps^{-m}}\abs{\partial^\alpha u_\eps(x)}\le\eps^{-N}\big),
\]
i.e., $\eps_0$ may be chosen independent of $m$.
\end{lemma}
\begin{proof}
Let $(\tilde u_\eps)_\eps\in \EMod(\GenR^d)$ be a representative of $u$. Then for each $m, k\in\N$, there exist $\eps_{m,k}\in (0,1)$ and $N_{m,k}\in\N$ such that
\[(\forall \alpha \text{ with }\abs\alpha\le k) (\forall\eps\le \eps_{m,k})\big(\sup_{\abs{x}\le \eps^{-m}}\abs{\partial^\alpha \tilde u_\eps(x)}\le\eps^{-N_{m,k}}\big).\]
We also introduce the notation $\eps_m:=\eps_{m,m}$.
We may suppose that $N_{m,k}$ is increasing in $m,k$, that $\eps_{m,k}$ is decreasing in $m,k$, and that $\lim_{m\to\infty}\eps_{m}=0$. Let $\chi\in\test(\R^d)$ with $\chi(x)=1$, if $\abs x\le 1/2$ and $\chi(x)=0$, if $\abs x\ge 1$. Let $\chi_\eps(x)=\chi(\eps^m x)$, if $\eps_{m+1}<\eps\le\eps_{m}$. Then $\chi_\eps(x)=1$, if $\abs x\le\eps^{-m}/2$ and $\chi_\eps(x)=0$, if $\abs x\ge \eps^{-m}$ (for $\eps_{m+1}<\eps\le\eps_m$). Further, $\sup_{\eps\in(0,1),x\in\R^d}\abs{\partial^\alpha \chi_\eps(x)}<\infty$, $\forall\alpha\in\N^d$. Let $u_\eps:= \tilde u_\eps\cdot \chi_\eps$. Then $u=[(u_\eps)_\eps]$.\\
Now let $\alpha\in\N^d$ be fixed. Choose $\eta\in(0,1)$ (depending on $\alpha$ only) with $\eta\le \eps_{\abs\alpha}$ and $2^{\abs\alpha}\sup_{\abs{\beta}\le\abs\alpha, \eps\in(0,1),  x\in\R^d} \abs{\partial^\beta \chi_\eps(x)}\le\eta^{-1}$. Then for each $m\in\N$ and $\eps\le\eta$,
\[
\sup_{\abs{x}\le \eps^{-m}}\abs{\partial^\alpha u_\eps(x)} \le\eta^{-1} \sup_{\abs\beta\le\abs\alpha, \abs{x}\le \eps^{-m}}\abs{\partial^\beta \tilde u_\eps(x)}\le\eps^{-1-N_{m,\abs\alpha}},
\]
as soon as $\eps\le\eps_{\max(m,\abs\alpha)}$.
If $\eps_{\max(m,\abs\alpha)}<\eps\le\eta$, then there exists $\abs\alpha\le k<\max(m,\abs\alpha)$ such that $\eps_{k+1}<\eps\le \eps_{k}$. Hence $\sup_{\abs{x}\le \eps^{-m}}\abs{\partial^\alpha u_\eps(x)} \le \sup_{\abs{x}\le \eps^{-k}}\abs{\partial^\alpha u_\eps(x)}\le\eps^{-1-N_{k,\abs\alpha}}$.
Summarizing, $\sup_{\abs{x}\le \eps^{-m}}\abs{\partial^\alpha u_\eps(x)}\le \eps^{-1-N_{\max(m,\abs\alpha),\abs\alpha}}$, as soon as $\eps\le\eta$.
\end{proof}

Our next aim is to characterize the subalgebras $\Gen_\tau(\R^d)$ and $\Gen_\Schwartz(\R^d)$ in $\Gen(\GenR^d)$ by the pointwise values of their elements (and their derivatives).
\begin{prop}\label{prop_E_tau_charac}
Let 
\begin{multline*}
\mathcal{V}=\{(u_\eps)_\eps\in\Cnt[\infty](\R^d)^{(0,1)}: (\forall\alpha\in\N^d) (\exists N\in\N) (\forall m\in\N\setminus\{0\})\\
(\exists \eps_0\in (0,1)) (\forall\eps\le\eps_0) (\sup_{\abs{x}\le\eps^{-m}}\abs{\partial^\alpha u_\eps(x)}\le\eps^{-mN})\}.
\end{multline*}
Then $\mathcal V = \E_\tau(\R^d)+\Null(\GenR^d)$
and $\Gen_\tau(\R^d)=\mathcal{V}/\Null(\GenR^d)$.
\end{prop}
\begin{proof}
Clearly, $\Null(\GenR^d)\subseteq\mathcal V$. Let $(u_\eps)_\eps\in\E_\tau(\R^d)$ and $\alpha\in\N^d$. Let $N\in\N$ such that $\sup_{x\in\R^d}(1+\abs{x})^{-N}\abs{\partial^\alpha u_\eps(x)}\le\eps^{-N}$, for small $\eps$. Then for each $m\in\N\setminus\{0\}$,
\[\sup_{\abs x\le\eps^{-m}}\abs{\partial^\alpha u_\eps(x)}\le\eps^{-N}\sup_{\abs x\le\eps^{-m}}(1+\abs{x})^{N}\le \eps^{-N}2^N\eps^{-mN}\le\eps^{-3Nm},\]
for small $\eps$. Hence $(u_\eps)_\eps\in\mathcal V$.
Since $\mathcal V$ is closed under $+$, $\E_\tau(\R^d)+\Null(\GenR^d)\subseteq \mathcal V$.\\
Conversely, let $(\tilde u_\eps)_\eps\in\mathcal V$. As in lemma \ref{lemma_eps_0_indep_of_m}, now with $N_k$ (independent of $m$) increasing in $k$, we find $(u_\eps)_\eps\in (\tilde u_\eps)_\eps + \Null(\GenR^d)$ such that for each $\alpha\in\N^d$, there exists $M=(\abs\alpha + 1)N_{\abs\alpha}+1\in\N$ and $\eta\in(0,1)$ such that
\[
(\forall m\in\N\setminus\{0\})
(\forall\eps\le\eta) \big(\sup_{\abs{x}\le\eps^{-m}}\abs{\partial^\alpha u_\eps(x)}\le\eps^{-mM}\big).
\]
Fix $\alpha\in\N^d$. Then for $\eps\le\eta$ and $m\in\N\setminus\{0\}$,
\[
\sup_{\eps^{-m} \le\abs x\le\eps^{-m-1}}(1+\abs{x})^{-M}\abs{\partial^\alpha u_\eps(x)}
\le \eps^{mM}\sup_{\abs x\le\eps^{-m-1}}\abs{\partial^\alpha u_\eps(x)}\le \eps^{-M}.
\]
Taking the supremum over all $m\in\N\setminus\{0\}$, we have for each $\eps\le\eta$ that
\[
\sup_{\abs x\ge \eps^{-1}}(1+\abs{x})^{-M}\abs{\partial^\alpha u_\eps(x)} \le \eps^{-M}.
\]
On the other hand, for $\eps\le\eta$,
\[
\sup_{\abs x\le \eps^{-1}}(1+\abs{x})^{-M}\abs{\partial^\alpha u_\eps(x)}
\le \sup_{\abs x\le \eps^{-1}}\abs{\partial^\alpha u_\eps(x)}\le\eps^{-M}.
\]
As $\alpha\in\N^d$ arbitrary, $(u_\eps)_\eps\in\E_\tau(\R^d)$.
\end{proof}

\begin{thm}\label{thm_E_tau_charac}
$\Gen_\tau(\R^d)=$
\[\{u\in\Gen(\GenR^d): (\forall\alpha\in\N^d) (\exists N\in\N) (\forall \tilde x\in\GenR^d\text{ with } \abs{\tilde x}\ge \caninf^{-1}) (\abs{\partial^\alpha u(\tilde x)}\le \abs{\tilde x}^N)\}.\]
\end{thm}
\begin{proof}
$\subseteq$: let $u\in\Gen_\tau(\R^d)$. Let $\alpha\in\N^d$. Since $u$ has a representative in $\E_\tau(\R^d)$, there exists $N\in\N$ such that for each $\tilde x\in\GenR^d$ with $\abs{\tilde x}\ge\caninf^{-1}$,  $\abs{\partial^\alpha u(\tilde x)} \le \caninf^{-N}(1+\abs{\tilde x})^N\le \caninf^{-N-1}\abs{\tilde x}^N\le\abs{\tilde x}^{2N+1}$.\\
$\supseteq$: if $u=[(u_\eps)_\eps]\in\Gen(\GenR^d)\setminus\Gen_\tau(\R^d)$, then, by proposition \ref{prop_E_tau_charac}, we find $\alpha\in\N^d$ such that for each $N\in\N$, we find $m\in\N\setminus\{0\}$, a decreasing sequence $(\eps_n)_n$ tending to $0$ and $x_{\eps_n}$ with $\abs{x_{\eps_n}}\le\eps_n^{-m}$ and $\abs{\partial^\alpha u_{\eps_n}(x_{\eps_n})}>\eps_n^{-mN}\ge\abs{x_{\eps_n}}^N$, $\forall n$. If $u\in\Gen(\GenR^d)$, there exists $M\in\N$ such that $\abs{\partial^\alpha u_\eps(x)}\le\eps^{-M}$, for $\abs{x}\le\eps^{-1}$ and for small $\eps$. Hence for $N>M$, $\abs{x_{\eps_n}}>\eps_n^{-1}$, for sufficiently large $n$. Let $x_\eps:=(\eps^{-1},0,\dots,0)$, if $\eps\notin\{\eps_n: n\in\N\}$. Then $\tilde x:=[(x_\eps)_\eps]\in\GenR^d$ with $\abs{\tilde x}\ge\caninf^{-1}$ and $\abs{\partial^\alpha u(\tilde x)}\nleq\abs{\tilde x}^{N-1}$.
\end{proof}

\begin{prop}\label{prop_E_Schwartz_charac}
Let 
\begin{multline*}
\mathcal{W}=\{(u_\eps)_\eps\in\EMod(\GenR^d): (\forall\alpha\in\N^d) (\forall k\in\N) (\exists N\in\N) (\forall m\in\N\setminus\{0\})\\
(\exists \eps_0\in (0,1)) (\forall\eps\le\eps_0) \big(\sup_{\eps^{-m}\le\abs{x}\le\eps^{-m-1}}\abs{\partial^\alpha u_\eps(x)}\le\eps^{mk-N}\big)\}.
\end{multline*}
Then $\mathcal W = \E_\Schwartz(\R^d)+\Null(\GenR^d)$
and $\Gen_\Schwartz(\R^d)=\mathcal{W}/\Null(\GenR^d)$.
\end{prop}
\begin{proof}
Clearly, $\Null(\GenR^d)\subseteq\mathcal W$. Let $(u_\eps)_\eps\in\E_\Schwartz(\R^d)$. Fix $\alpha\in\N^d$ and $k\in \N$. Then there exists $N\in\N$ such that $\sup_{x\in\R^d}\abs{x}^{k}\abs{\partial^\alpha u_\eps(x)}\le\eps^{-N}$, for small $\eps$. Thus for each $m\in\N\setminus\{0\}$,
\[
\sup_{\eps^{-m}\le\abs x\le\eps^{-m-1}}\abs{\partial^\alpha u_\eps(x)} \le\eps^{-N}\sup_{\eps^{-m}\le\abs x\le\eps^{-m-1}}\abs{x}^{-k}
\le \eps^{mk-N},
\]
for small $\eps$. Hence $(u_\eps)_\eps\in\mathcal W$.
Since $\mathcal W$ is closed under $+$, $\E_\Schwartz(\R^d)+\Null(\GenR^d)\subseteq \mathcal W$.\\
Conversely, let $(\tilde u_\eps)_\eps\in\mathcal W$. As in lemma \ref{lemma_eps_0_indep_of_m}, now with $\eps_{\abs\alpha,k,m}$ and $\eps_m:=\eps_{m,m,m}$, we find $(u_\eps)_\eps\in (\tilde u_\eps)_\eps + \Null(\GenR^d)$ such that for each $\alpha\in\N^d$ and $k\in\N$, there exists $N\in\N$ and $\eta\in(0,1)$ (choose now $\eta\le\eps_{\max(\abs\alpha,k)}$) such that
\[
(\forall m\in\N\setminus\{0\})
(\forall\eps\le\eta) \big(\sup_{\eps^{-m}\le\abs{x}\le\eps^{-m-1}}\abs{\partial^\alpha u_\eps(x)}\le\eps^{mk-N}\big).
\]
Fix $\alpha\in\N^d$ and $k\in\N$. Then for $\eps\le\eta$, $\beta\in\N^d$ with $\abs\beta=k$ and $m\in\N\setminus\{0\}$,
\[
\sup_{\eps^{-m} \le\abs x\le\eps^{-m-1}}\abs{x^\beta\partial^\alpha u_\eps(x)}
\le \eps^{-(m+1)\abs\beta}\sup_{\eps^{-m}\le\abs x\le\eps^{-m-1}}\abs{\partial^\alpha u_\eps(x)}\le \eps^{-N-\abs\beta}.
\]
Taking the supremum over all $m\in\N\setminus\{0\}$, we have for each $\eps\le\eta$ that
\[
\sup_{\abs x\ge \inv\eps}\abs{x^\beta\partial^\alpha u_\eps(x)} \le \eps^{-N-\abs\beta}.
\]
On the other hand, since $(u_\eps)_\eps\in\EMod(\GenR^d)$, there exists $M\in\N$ such that
\[
\sup_{\abs x\le\inv\eps}\abs{x^\beta\partial^\alpha u_\eps(x)}
\le \eps^{-\abs\beta}\sup_{\abs x\le \inv\eps}\abs{\partial^\alpha u_\eps(x)}\le\eps^{-M-\abs\beta},
\]
for sufficiently small $\eps$. As $\alpha,\beta\in\N^d$ arbitrary, $(u_\eps)_\eps\in\E_\Schwartz(\R^d)$.
\end{proof}

\begin{thm}\label{thm_E_Schwartz_charac}
$\Gen_\Schwartz(\R^d)=\{u\in\Gen(\GenR^d):$
\[(\forall\alpha\in\N^d) (\forall k\in\N) (\exists C\in\GenR) (\forall \tilde x\in\GenR^d\text{ with } \abs{\tilde x}\ge \inv\caninf) (\abs{\partial^\alpha u(\tilde x)}\le C\abs{\tilde x}^{-k})\}.\]
\end{thm}
\begin{proof}
$\subseteq$: Let $u\in\Gen_\Schwartz(\R^d)$. Let $\alpha\in\N^d$ and $k\in\N$. Since $u$ has a representative in $\E_\Schwartz(\R^d)$, we find $N\in\N$ such that for each $\beta\in\N^d$ with $\abs\beta\le k$ and $\tilde x\in\GenR^d$, $\abs{\tilde x^\beta \partial^\alpha u(\tilde x)}\le \caninf^{-N}$. Hence $\abs{\tilde x}^k\abs{\partial^\alpha u(\tilde x)}\le(\sum_{j=1}^d\abs{x_j})^k\abs{\partial^\alpha u(\tilde x)}\le d^k\caninf^{-N}$. If $\abs{\tilde x}\ge \inv\caninf$, this implies $\abs{\partial^\alpha u(\tilde x)}\le C\abs{\tilde x}^{-k}$, with $C=d^k\caninf^{-N}\in\GenR$.\\
$\supseteq$: if $u=[(u_\eps)_\eps]\in\Gen(\GenR^d)\setminus\Gen_\Schwartz(\R^d)$, then, by proposition \ref{prop_E_Schwartz_charac}, we find $\alpha\in\N^d$ and $k\in\N$ such that for each $N\in\N$, we find $m\in\N\setminus\{0\}$, a decreasing sequence $(\eps_n)_{n\in\N}$ tending to $0$ and $x_{\eps_n}$ with $\eps_n^{-m}\le\abs{x_{\eps_n}}\le\eps_n^{-m-1}$ and $\abs{\partial^\alpha u_{\eps_n}(x_{\eps_n})}>\eps_n^{mk-N}\ge \abs{x_{\eps_n}}^{-k}\eps_n^{-N}$.
Let $x_\eps:=(\inv\eps,0,\dots,0)$, if $\eps\notin\{\eps_n: n\in\N\}$. Then $\tilde x:=[(x_\eps)_\eps]\in\GenR^d$ with $\abs{\tilde x}\ge \inv\caninf$ and $\abs{\partial^\alpha u(\tilde x)}\nleq\abs{\tilde x}^{-k}\caninf^{-N+1}$. As $N\in\N$ is arbitrary, we conclude that for each $C\in\GenR$, there exists $\tilde x\in\GenR^d$ with $\abs{\tilde x}\ge \inv\caninf$ such that $\abs{\partial^\alpha u(\tilde x)}\nleq C\abs{\tilde x}^{-k}$.
\end{proof}

\begin{df}
A generalized point $\tilde x\in\GenR^d$ is said to be \defstyle{of slow scale} if $\sharpnorm{\tilde x}\le 1$, or equivalently, if for each $a\in\R^+$, 
$\abs{x_\eps}\le\eps^{-a}$, for small $\eps$. Similarly, a generalized point $\tilde x\in\GenR^d$ is said to be \defstyle{of fast scale} if $\abs{\tilde x}$ is invertible in $\GenR$ and $\sharpnorm{1/\tilde x} < 1$, or equivalently, if there exists $a\in\R^+$ such that $\abs{x_\eps}\ge\eps^{-a}$, for small $\eps$. We denote the set of slow, resp.\ fast, scale points of $\GenR^d$ by $\GenR^d_{ss}$, resp.\ $\GenR^d_{fs}$.
\end{df}
\begin{thm}\label{thm_slowscale_support}
For $u\in\Gen(\GenR^d)$, the following are equivalent:
\begin{enumerate}
\item $u(\tilde x)=0$, $\forall \tilde x\in\GenR^d_{fs}$
\item $\partial^\alpha u(\tilde x)=0$, $\forall \tilde x\in\GenR^d_{fs}$, $\forall \alpha\in\N^d$
\item $(\exists \tilde a\in\GenR_{ss})$ $(\forall \tilde x\in\GenR^d$ with $\abs{\tilde x}\ge\tilde a)$ $(u(\tilde x)=0)$
\item $u$ has a representative $(u_\eps)_\eps$ such that
\[(\forall a\in\R^+) (\exists\eta\in(0,1)) (\forall\eps\le\eta) (\forall x\in\R^d\text{ with }\abs x\ge \eps^{-a}) (u_\eps(x) = 0)
\]
\item $u$ has a representative $(u_\eps)_\eps$ such that
\[(\forall m\in\N\setminus\{0\})
(\exists\eps_m\in(0,1))(\forall\eps\le\eps_m)\big(\sup_{\abs x\ge \eps^{-1/m}} \abs{u_\eps(x)}\le\eps^{m}\big)\]
\item $u$ has a representative $(u_\eps)_\eps$ such that
\[(\forall\alpha\in\N^d) (\forall m\in\N\setminus\{0\})
(\exists\eta\in(0,1))(\forall\eps\le\eta)\big(\sup_{\abs x\ge \eps^{-1/m}} \abs{\partial^\alpha u_\eps(x)}\le\eps^{m}\big)\]
\item $u$ has a representative $(u_\eps)_\eps$ such that
\[
(\exists N\in\N)(\forall\beta\in\N^d)\\
(\exists\eta\in(0,1))(\forall\eps\le\eta)(\sup_{x\in\R^d} \abs{x^\beta u_\eps(x)}\le\eps^{-N})\]
\item $u$ has a representative $(u_\eps)_\eps$ such that
\[
(\forall\alpha\in\N^d) (\exists N\in\N)(\forall\beta\in\N^d)\\
(\exists\eta\in(0,1))(\forall\eps\le\eta)(\sup_{x\in\R^d} \abs{x^\beta\partial^\alpha u_\eps(x)}\le\eps^{-N})\]
(in particular, $u$ has a representative in $\mathcal{E}_\Schwartz(\R^d)$).
\end{enumerate}
We call \defstyle{slow scale supported} those $u\in\Gen(\GenR^d)$ satisfying one of these equivalent conditions, and denote the set of all slow scale supported $u\in\Gen(\GenR^d)$ by $\Gen_{ss}(\GenR^d)$.\\
Then $\Gen_\Schwartz^\infty(\R^d) \subsetneqq\Gen_{ss}(\GenR^d) \subsetneqq\Gen_\Schwartz(\R^d)$ is a differential ideal of $\Gen(\GenR^d)$.
\end{thm}
\begin{proof}
$(1)\implies(2)$: Let $\tilde x_0\in\GenR^d_{fs}$. Then $v(x):=u(x+\tilde x_0)\in\Gen(\GenR^d)$ and $v(\tilde x)=0$ for each $\tilde x\in\GenR^d_c$. By the pointwise value theorem in $\Gen(\R^d)$, $\restr{v}{\R^d}=0$ in $\Gen(\R^d)$. Hence also $\partial^\alpha u(\tilde x_0)=\partial^\alpha v(0)=0$, $\forall\alpha\in\N^d$.\\
$(2)\implies(6)$: By contraposition, $(2)$ implies that for a representative $(u_\eps)_\eps$ of $u$,
\[(\forall\alpha\in\N^d) (\forall k\in\N) (\forall m\in\N\setminus\{0\}) (\exists \eps_0\in (0,1)) (\forall \eps\le\eps_0) (\sup_{\eps^{-m}\le\abs x\le \eps^{-m-1}}\abs{\partial^\alpha u_\eps(x)}\le\eps^k).\]
With the notations of proposition \ref{prop_E_Schwartz_charac}, this implies that $(u_\eps)_\eps\in\mathcal W$, hence $u\in\Gen_\Schwartz(\R^d)$. So we may assume that $(u_\eps)_\eps\in\mathcal{E}_\Schwartz(\R^d)$.
Now suppose that $(6)$ does not hold, then
we find $\alpha\in\N^d$, $M\in\N\setminus\{0\}$, a decreasing sequence $(\eps_n)_{n\in\N}$ tending to $0$ and $x_{\eps_n}\in\R^d$ with $\abs{x_{\eps_n}}\ge\eps_n^{-1/M}$ such that $\abs{\partial^\alpha u_{\eps_n}(x_{\eps_n})}>\eps_n^M$, $\forall n$. As $(u_\eps)_\eps\in\E_\Schwartz(\R^d)$, there exists $N\in\N$ such that
$\sup_{x\in\R^d}\abs{x} \abs{\partial^\alpha u_{\eps}(x)}\le \eps^{-N}$, for small $\eps$. Hence $\abs{x_{\eps_n}}\le\eps_n^{-N} \abs{\partial^\alpha u_{\eps_n}(x_{\eps_n})}^{-1} < \eps_n^{-N-M}$, for small $n$. Let $x_\eps\in\R^d$ with $\abs {x_\eps} = \eps^{-1/M}$, if $\eps\notin\{\eps_n: n\in\N\}$. Then $\tilde x:=[(x_\eps)_\eps]\in\GenR^d_{fs}$ but $\partial^\alpha u(\tilde x)\ne 0$, contradicting $(2)$.\\
$(6)\implies(5)$: trivial.\\
$(5)\implies(3)$: For each $m\in\N\setminus\{0\}$, let $\eps_m$ as in the statement of $(5)$. We may assume that $(\eps_m)_{m\in\N}$ is decreasing and tends to $0$. Let $a_\eps=\eps^{-1/\sqrt m}$, for $\eps_{m+1}<\eps\le\eps_m$. As $a_\eps\ge \eps^{-1/m}$, $\abs{u_\eps(x)}\le\eps^m$, for each $x\in\R^d$ with $\abs x\ge a_\eps$ and $\eps_{m+1}<\eps\le\eps_m$. Thus $\tilde a:=[(a_\eps)_\eps]\in\GenR$ satisfies the statement of $(3)$.\\
$(3)\implies(1)$: If $\tilde x\in\GenR^d_{fs}$, then there exists $b\in\R^+$ such that $\abs{\tilde x}\ge\caninf^{-b}\ge \tilde a$.\\
$(3)\implies(4)$: Let $\tilde a=[(a_{\eps})_\eps]$ as in the statement of (3). We may suppose that $a_{\eps}\ge 1$, $\forall\eps$.
Let $(\chi_\eps)_\eps\in\EMod(\GenR^d)$ with $\chi_\eps(x)=1$, if $\abs x\le a_{\eps}$ and $\chi_\eps(x)=0$, if $\abs x\ge 2 a_{\eps}$. Let $u=[(\tilde u_\eps)_\eps]$ and define $u_\eps:=\chi_\eps\cdot\tilde u_\eps$. Then $(u_\eps-\tilde u_\eps)_\eps \in\Null(\GenR^d)$ since the point-values in all points of $\GenR^d$ are $0$ (to be precise, one uses an interleaving argument as in the proof of lemma \ref{lemma_support}). Hence $u=[(u_\eps)_\eps]$.\\
$(4)\implies(8)$: Fix $\alpha\in\N^d$. As $(u_\eps)_\eps\in\EMod(\GenR^d)$, there exists $N\in\N$ such that, for small $\eps$, $\sup_{\abs x\le\eps^{-1}}\abs{\partial^\alpha u_\eps(x)}\le\eps^{-N}$. Then for $\beta\in\N^d$,
\[
\sup_{x\in\R^d}\abs{x^\beta\partial^\alpha u_\eps(x)}
\le \sup_{\abs x\le \eps^{-1/(\abs\beta+1)}}\abs{x^\beta\partial^\alpha u_\eps(x)}\le \eps^{-\frac{\abs\beta}{\abs\beta + 1}}\eps^{-N}\le\eps^{-N-1},
\]
for small $\eps$.\\
$(8)\implies(7)$: trivial.\\
$(7)\implies(5)$: There exists $N\in\N$ such that for each $k\in\N$, $\sup_{x\in\R^d}\abs{x}^k\abs{u_\eps(x)}\le\eps^{-N}$, for small $\eps$. Hence for each $k\in\N$ and $m\in\N\setminus\{0\}$,
\[
\sup_{\abs x\ge\eps^{-1/m}}\abs{u_\eps(x)}
\le \eps^{-N}\sup_{\abs x\ge\eps^{-1/m}}\abs{x}^{-k}\le \eps^{k/m-N},
\]
for small $\eps$.
\end{proof}

Since $\Gen_{ss}(\GenR^d)\subseteq\Gen_\Schwartz(\R^d)$, the Fourier transform $\widehat{\ }$: $\Gen_{ss}(\GenR^d)\to\Gen_\Schwartz(\R^d)$ is well-defined and $\widehat u = [(\widehat u_\eps)_\eps]$, if $(u_\eps)_\eps\in\E_\Schwartz(\R^d)$ \cite{GarettoPhD}.

\begin{thm}\label{thm_slowscale_spectrum}
The Fourier transform $\widehat\Gen_{ss}(\GenR^d)$ of $\Gen_{ss}(\GenR^d)$, consisting of those elements in $\Gen(\GenR^d)$ with \defstyle{slow scale spectrum}, is given by those elements in $\Gen_\Schwartz(\R^d)$ with a representative $(u_\eps)_\eps\in\E_\Schwartz(\R^d)$ satisfying
\begin{equation}\label{eqn_G_ss_Fourier}
(\forall\alpha\in\N^d) (\exists N\in\N)(\forall\beta\in\N^d)\\
(\exists\eta\in(0,1))(\forall\eps\le\eta)(\sup_{x\in\R^d} \abs{x^\alpha\partial^\beta u_\eps(x)}\le\eps^{-N}).
\end{equation}
Further, $\widehat\Gen_{ss}(\R^d)$ is a differential subalgebra of $\Gen_\Schwartz(\R^d)$ and
\begin{align*}
\Gen_\Schwartz^\infty(\R^d) &= \Gen_{ss}(\GenR^d)\cap\widehat\Gen_{ss}(\GenR^d)
=\{u\in\Gen(\GenR^d): u(\tilde x)=\widehat u(\tilde x)=0, \quad\forall \tilde x\in\GenR^d_{fs}\}.
\end{align*}
\end{thm}
\begin{proof}
Let $(u_\eps)_\eps\in\E_\Schwartz(\R^d)$ satisfy equation (\ref{eqn_G_ss_Fourier}). Using characterization (8) of theorem \ref{thm_slowscale_support}, $\widehat u\in\Gen_{ss}(\GenR^d)$, since
\begin{multline*}
\sup_{\xi\in\R^d}\abs{\xi^\beta \partial^\alpha \widehat u_\eps(\xi)}
=\sup_{\xi\in\R^d}\abs{\widehat{\partial^\beta (x^\alpha u_\eps)}(\xi)}
\le \sum_{\gamma\le\beta} \int_{\R^d}\abs{\partial^{\beta-\gamma}x^\alpha}\, \abs{\partial^{\gamma}u_\eps(x)}\,dx\\
\le \eps^{-1} \sup_{\gamma\le\beta,x\in\R^d} (1+\abs{x})^{d+\abs\alpha+1} \abs{\partial^\gamma u_\eps(x)},
\end{multline*}
for small $\eps$. By Fourier inversion, $u\in\widehat\Gen_{ss}(\GenR^d)$. Conversely, by the same inequality, if $u\in\Gen_{ss}(\GenR^d)$, then there exists a representative $(u_\eps)_\eps\in\E_\Schwartz(\R^d)$ of $u$ such that $(\widehat u_\eps)_\eps\in\E_\Schwartz(\R^d)$ satisfies equation (1). The result follows again by Fourier inversion.

Finally, let $u\in\Gen_{ss}(\GenR^d)\cap\widehat\Gen_{ss}(\GenR^d)$. Let $(u_\eps)_\eps\in\E_\Schwartz(\R^d)$ be a representative of $u$ satisfying condition (4) of theorem \ref{thm_slowscale_support}.
Since equation (\ref{eqn_G_ss_Fourier}) holds independent of the representative in $\E_\Schwartz(\R^d)$, there exists $N\in\N$ such that for each $\beta\in\N^d$ and for small $\eps$, $\sup_{x\in\R^d}\abs{\partial^\beta u_\eps(x)}\le\eps^{-N}$. Then, for any $\alpha$, $\beta\in\N^d$,
\[
\sup_{x\in\R^d} \abs{x^\alpha \partial^\beta u_\eps(x)}
\le \sup_{\abs x\le \eps^{-1/(\abs\alpha + 1)}} \abs{x^\alpha \partial^\beta u_\eps(x)}
\le \eps^{-N-1},
\]
for small $\eps$.
\end{proof}

\section{Equalities in the sense of tempered generalized distributions}
\label{sec_gtd}
Let $u,v\in\Gen_\Schwartz(\R^d)$. In this section, we investigate when $u$ is equal to $v$ in the sense of generalized tempered distributions \cite{Colombeau85}, i.e., when $\int_{\R^d} u\phi = \int_{\R^d} v\phi$, for each $\phi\in\Schwartz(\R^d)$.

\begin{lemma}\label{lemma_support}
Let $u\in\Gen_{ss}(\GenR^d)$ and $v\in\Gen(\GenR^d)$ with $v(\tilde x)=0$, $\forall \tilde x\in\GenR^d_{ss}$, then $uv=0$ in $\Gen(\GenR^d)$.
\end{lemma}
\begin{proof}
By theorem \ref{thm_slowscale_support}, there exists $\tilde x_0\in\GenR^d_{ss}$ such that $u(\tilde y)=0$, for each $\tilde y\in\GenR^d$ with $\abs{\tilde y}\ge \abs{\tilde x_0}$. Let $\tilde x\in\GenR^d$. Fix representatives $(x_\eps)_\eps$ of $\tilde x$ and $(x_{0,\eps})_\eps$ of $\tilde x_0$, and let $S=\{\eps\in (0,1): \abs{x_\eps}\le \abs{x_{0,\eps}}\}$. Then $\abs{\tilde x e_{\co S}+ \tilde x_0 e_S}\ge\abs{\tilde x_0}$ and $\sharpnorm{\tilde x e_S}\le \sharpnorm{\tilde x_0}\le 1$ ($e_S\in\GenR$ is the element with the characteristic function on $S$ as a representative). Hence $u(\tilde x e_{\co S}+ \tilde x_0 e_S)=u(\tilde x) e_{\co S} + u(\tilde x_0) e_S=0$, $v(\tilde x e_S)=v(\tilde x)e_S=0$, so $u(\tilde x)v(\tilde x)=u(\tilde x)(v(\tilde x)e_S)+ (u(\tilde x)e_{\co S})v(\tilde x)= 0$. As $\tilde x\in\GenR^d$ arbitrary, $uv=0$ in $\Gen(\GenR^d)$.
\end{proof}

\begin{lemma}\label{lemma_product_insertion}
Let $u\in\Gen_{ss}(\GenR^d)$. Then there exists $\chi\in\Gen_\Schwartz^\infty(\R^d)$ such that $u=u\chi$.
\end{lemma}
\begin{proof}
By theorem \ref{thm_slowscale_support}, there exists $\tilde a\in\GenR_{ss}$ with $\tilde a\ge 1$ such that $u(\tilde x)=0$, for each $\tilde x\in\GenR^d$ with $\abs{\tilde x}\ge \tilde a$. Let $\tilde a = [(a_{\eps})_\eps]$ with $a_{\eps}\ge 1$, $\forall \eps$ and let $\phi\in\test(\R^d)$ with $\phi(x)=1$, if $\abs x\le 1$ and $\phi(x)=0$, if $\abs x\ge 2$. Let $\chi_\eps(x)=\phi(x/a_{\eps})$. It is easy to see that $\chi=[(\chi_\eps)_\eps]\in\Gen_{\Schwartz}^\infty(\R^d)$.
Further, $u(\tilde x)=u(\tilde x)\chi(\tilde x)$, for each $\tilde x\in\GenR^d$, as can be seen analogously to lemma \ref{lemma_support}. Hence $u=u\chi$ in $\Gen(\GenR^d)$ (and hence also in $\Gen_\Schwartz(\R^d)$).
\end{proof}

\begin{prop}\label{prop_distr_eq_zero_ss}
Let $u\in\Gen_\Schwartz(\R^d)$. Then the following are equivalent:
\begin{enumerate}
\item $\int_{\R^d}u\phi=0$, $\forall\phi\in\widehat\Gen_{ss}(\GenR^d)$
\item $u\conv\phi=0$, $\forall\phi\in\widehat\Gen_{ss}(\GenR^d)$
\item $\widehat u\phi=0$, $\forall\phi\in\Gen_{ss}(\GenR^d)$
\item $\widehat u(\tilde \xi)=0$, $\forall \tilde\xi\in\GenR^d_{ss}$
\item $\int_{\R^d}\widehat u\phi=0$, $\forall\phi\in\Gen_{ss}(\GenR^d)$.
\end{enumerate}
\end{prop}
\begin{proof}
$(1)\Rightarrow(2)$: if $\phi\in\widehat\Gen_{ss}(\GenR^d)$, then also $x\mapsto\phi(\tilde a-x)\in\widehat\Gen_{ss}(\GenR^d)$, for each $\tilde a\in\GenR^d$.\\
$(2)\Rightarrow(3)$: $\widehat u\cdot \widehat\phi = \widehat{u\conv\phi}=\widehat 0 = 0$, $\forall\phi\in\widehat\Gen_{ss}(\GenR^d)$.\\
$(3)\Rightarrow(4)$: let $\tilde\xi\in\GenR^d_{ss}$. Let $\phi\in\Schwartz(\R^d)$ with $\phi(0)=1$. Then $\psi(x):=\phi(x-\tilde\xi)\in\Gen_{ss}(\GenR^d)$ and $\widehat u(\tilde\xi)=\widehat u(\tilde\xi) \psi(\tilde\xi)=0$.\\
$(4)\Rightarrow(5)$: by lemma \ref{lemma_support}, $\widehat u\phi=0$, $\forall\phi\in\Gen_{ss}(\GenR^d)$.\\
$(5)\Rightarrow(1)$: by Parseval's formula, $\int_{\R^d}\widehat u \phi=\int_{\R^d}u\widehat \phi$, $\forall u,\phi\in\Gen_\Schwartz(\R^d)$.
\end{proof}

\begin{thm}\label{thm_gds_sufficient}
Let $u\in\Gen_\Schwartz(\R^d)$. Let $u = v + w$, where $v, w\in\Gen_\Schwartz(\R^d)$ and $v(\tilde x)=\widehat w(\tilde x) = 0$, $\forall\tilde x\in\GenR^d_{ss}$. Then $\int_{\R^d} u\phi=0,$ $\forall\phi\in\Schwartz(\R^d)$. 
\end{thm}
\begin{proof}
By proposition \ref{prop_distr_eq_zero_ss}, $\int_{\R^d} w\phi=0$, $\forall\phi\in\Schwartz(\R^d)\subseteq\widehat\Gen_{ss}(\GenR^d)$.\\
If $\phi\in\Schwartz(\R^d)\subseteq \Gen_{ss}(\GenR^d)$, then $v\phi=0$ in $\Gen_\Schwartz(\R^d)$ by lemma \ref{lemma_support}. Hence $\int_{\R^d} v\phi=0$.
\end{proof}

\begin{conj}
There exists $u\in\Gen_\Schwartz(\R^d)$ such that $\int_{\R^d} u\phi=0$, $\forall\phi\in\Schwartz(\R^d)$, yet $u\ne v + w$, where $v,w\in\Gen_\Schwartz(\R^d)$ and $v(\tilde x)=\widehat w(\tilde x) = 0$, $\forall\tilde x\in\GenR^d_{ss}$.
\end{conj}

If $u\in\Gen_{ss}(\GenR^d)$ or $u\in\widehat\Gen_{ss}(\GenR^d)$, then the converse of theorem \ref{thm_gds_sufficient} holds.
\begin{thm}\label{thm_gds_characterization}\leavevmode
\begin{enumerate}
\item If $u\in\Gen_{ss}(\GenR^d)$, then
\[\int_{\R^d} u\phi=0,\quad\forall\phi\in\Schwartz(\R^d) \iff \widehat u(\tilde \xi)=0, \quad \forall \tilde \xi\in\GenR^d_{ss}.\]
\item If $u\in\widehat\Gen_{ss}(\GenR^d)$, then
\[\int_{\R^d} u\phi=0, \quad\forall\phi\in\Schwartz(\R^d) \iff u(\tilde x)=0, \quad \forall \tilde x \in\GenR^d_{ss}.\]
\end{enumerate}
\end{thm}
\begin{proof}
(1) $\Rightarrow$: as $\Schwartz(\R^d)$ is a barreled topological vector space, also $\int_{\R^d}u\phi=0$, $\forall\phi\in\Gen_\Schwartz^\infty(\R^d)$ by \cite[Thm.~3.3]{HV_weak_homog} (this follows as in \cite[Thm.~3.5]{HV_weak_homog}).
Let $\psi\in\widehat\Gen_{ss}(\GenR^d)$. By lemma \ref{lemma_product_insertion}, $u=u\chi$ with $\chi\in\Gen_\Schwartz^\infty(\R^d)$. Hence $\int_{\R^d}u\psi=\int_{\R^d}u(\chi\psi)=0$, since $\chi\psi\in\Gen_{ss}(\GenR^d)\cap\widehat\Gen_{ss}(\GenR^d)=\Gen_\Schwartz^\infty(\R^d)$. The conclusion follows by proposition \ref{prop_distr_eq_zero_ss}.\\
(2) $\Rightarrow$: by Parseval's formula, $\int_{\R^d} \widehat u\phi =\int_{\R^d} u \widehat\phi=0$, $\forall\phi\in\Schwartz(\R^d)$. Since $\widehat u\in\Gen_{ss}(\GenR^d)$, the conclusion follows by part (1) and Fourier inversion.\\
(1), (2) $\Leftarrow$: by theorem \ref{thm_gds_sufficient}.
\end{proof}
Since $\Gen_c(\R^d)\subseteq\Gen_{ss}(\GenR^d)$, theorem \ref{thm_gds_characterization} can be viewed as a generalization of \cite[Thm.~4.2]{HV_weak_homog}.

\section{Pointwise $\Gen^\infty$-regularity}\label{sec_regularity}
We recall the various definitions of a generalized function regular at a non-generalized point and at a compactly supported generalized point as introduced in \cite{DPS,OPS} (in \cite{OPS}, other than compactly supported generalized points are considered as well).
\begin{df}\cite{DPS,OPS}
Let $u=[(u_\eps)_\eps]\in\Gen(\Omega)$. Let $x_0\in\Omega$ and $\tilde x_0\in\widetilde\Omega_c$.
\begin{enumerate}
\item $u\in\Gen^\infty_{x_0}$ iff
$(\exists V\in \mathcal V_{x_0}) (\exists N\in\N) (\forall\alpha\in\N^d) \big(\sup\limits_{x\in V}\abs{\partial^\alpha u_\eps(x)}\le\eps^{-N}$, for small $\eps\big)$.
\item $u\in\dot\Gen^\infty_{x_0}$ iff
$(\exists N\in\N) (\forall\alpha\in\N^d) (\exists V\in \mathcal V_{x_0}) \big(\sup\limits_{x\in V}\abs{\partial^\alpha u_\eps(x)}\le\eps^{-N}$, for small $\eps\big)$.
\item $u\in\widetilde\Gen^\infty_{\widetilde x_0}$ iff
$(\exists V\in \widetilde{\mathcal V}_{\widetilde x_0}) (\exists N\in\N) (\forall\alpha\in\N^d) (\forall \tilde x\in V) (\abs{\partial^\alpha u(\tilde x)}\le\caninf^{-N})$.
\item $u\in\Gentdinfty_{\widetilde x_0}$ iff
$(\exists N\in\N) (\forall\alpha\in\N^d) (\abs{\partial^\alpha u(\tilde x_0)}\le\caninf^{-N})$.
\end{enumerate}
Further, $\dot\Gen^\infty(\Omega):=\bigcap_{x\in\Omega}\dot\Gen^\infty_{x}$. The equality $\Gen^\infty(\Omega)=\bigcap_{x\in\Omega}\Gen^\infty_{x}$ holds.
\end{df}

We extend the definition of $\Gentdinfty$ at a point to arbitrary subsets of $\widetilde\Omega_c$ and indicate how the previous definitions can be related to each other by means of this definition. For this purpose, we extend (by similar arguments) the result given in \cite{OPS} that $u\in\Gen^\infty(\Omega)$ iff $u\in\Gentdinfty_{\widetilde x}$ for each $\tilde x\in\widetilde\Omega_c$.

\begin{df}
Let $\emptyset\ne A\subseteq\widetilde\Omega_c$ and $u\in\Gen(\Omega)$. We say that $u\in\Gentdinfty(A)$ iff $u\in\Gentdinfty_{\widetilde x}$ for each $\tilde x\in A$.
\end{df}

We recall \cite{OV_internal} that $A\subseteq\GenR^d$ is called \defstyle{internal} iff there exists a net $(A_\eps)_\eps$ of subsets of $\R^d$ such that $A$ is the set of those $\tilde x\in \GenR^d$ with a representative $(x_\eps)_\eps$ such that $x_\eps\in A_\eps$, for small $\eps$. In this case, we write $A=[(A_\eps)_\eps]$. In particular, for $A\subseteq\R^d$, we denote by $\widetilde A$ the internal set $[(A)_\eps]$.

\begin{prop}\label{prop_geninftyinternal}
Let $(A_n)_{n\in\N}$ be a decreasing sequence of internal subsets of $\widetilde\Omega$. Let $A_n=[(A^{n}_{\eps})_\eps]$, for each $n\in\N$. Let $B = \bigcap_{n\in\N} A_n$. Let $B_c=B\cap\widetilde\Omega_c\ne\emptyset$ and $u=[(u_\eps)_\eps]\in\Gen(\Omega)$.
Then the following are equivalent.
\begin{enumerate}
\item $u\in\Gentdinfty(B_c)$
\item $(\forall K\csub \Omega) (\exists N\in\N) (\forall\alpha\in\N^d) (\exists m\in\N) (\forall \tilde x\in \widetilde K\cap A_m) (\abs{\partial^\alpha u(\tilde x)}\le\caninf^{-N})$
\item $(\forall K\csub\Omega) (\exists N\in\N) (\forall\alpha\in\N^d) (\exists m\in\N) \big(\sup\limits_{x\in K\cap A^{m}_\eps}\abs{\partial^\alpha u_\eps(x)}\le\eps^{-N}, \text{ for small }\eps\big)$.
\end{enumerate}
\end{prop}
\begin{proof}
$(1)\Rightarrow(3)$: Supposing the conclusion is not true, we find $K\csub\Omega$, $\alpha_n\in\N^d$ (for each $n\in\N$), $\eps_{n,m}\in (0,1/m)$ (for each $n,m\in\N$) (by enumerating the countable family $(\eps_{n,m})_{n,m}$, we can successively choose the $\eps_{n,m}$ in such a way that they are all different) and $x_{\eps_{n,m}}\in K\cap A^{n+m}_{\eps_{n,m}}$ with $\abs{\partial^{\alpha_n}u_{\eps_{n,m}}(x_{\eps_{n,m}})} > \eps_{n,m}^{-n}$, $\forall n,m\in\N$. Let $\tilde y=[(y_\eps)_\eps]\in B_c$. Then there exists $L\csub\Omega$ such that $K\subseteq L$ and $\tilde y\in [(L\cap A_\eps^k)_\eps]$, for each $k\in\N$. For $k\le n+m$, $[(L\cap A_\eps^{n+m})_\eps]\subseteq A_k$, so by \cite[Prop.~2.9]{OV_internal}, $\sup_{x\in L\cap A_\eps^{n+m}}d(x,A_\eps^k)\le\eps^{m+n}$, for small $\eps$. Hence we can ensure that $d(x_{\eps_{n,m}},A_{\eps_{n,m}}^k)\le\eps_{n,m}^{m+n}$, for each $n,m\in\N$ and $k\le n+m$. Let $x_\eps=y_\eps$, if $\eps\notin\{\eps_{n,m}: n,m\in\N\}$.
Let $k,l\in\N$. Then $d(x_{\eps_{n,m}}, A_{\eps_{n,m}}^k)\le\eps_{n,m}^l$, except for finitely many $(n,m)\in\N^2$. Since also $\tilde y\in A_k$, $(d(x_\eps, A_\eps^k))_\eps\in\Null_\R$ by \cite[Prop.~2.1]{OV_internal}. By the same proposition, $\tilde x=[(x_\eps)_\eps]\in B_c$. By hypothesis, there exists $N\in\N$ such that for each $\alpha\in\N^d$, $\abs{\partial^\alpha u(\tilde x)}\le\caninf^{-N}$. This contradicts the fact that for a fixed $n>N$, $\lim_{m\to\infty}\eps_{n,m}=0$ and $\abs{\partial^{\alpha_n} u_{\eps_{n,m}}(x_{\eps_{n,m}})}> \eps_{n,m}^{-n}$, $\forall m\in\N$.\\
$(3)\Rightarrow(2)$: let $K$ be contained in the interior of $L\csub\Omega$. Then for each $m\in\N$ and $\tilde x\in \widetilde K\cap A_m$, there exists a representative $(x_\eps)_\eps$ of $\tilde x$ with $x_\eps\in L\cap A^m_\eps$, for small $\eps$.\\
$(2)\Rightarrow(1)$: $\tilde x\in B_c$ iff $\tilde x\in B\cap \widetilde K$, for some $K\csub\Omega$.
\end{proof}

\begin{df}
Let $\tilde x,\tilde y\in \widetilde\Omega_c$. We say that $\tilde x$ is infinitely close to $\tilde y$ (notation: $\tilde x\approx \tilde y$) iff $\abs{\tilde x-\tilde y}<1/n$, for each $n\in\N$ (in Colombeau theory, it is also said that $x-y$ is associated with $0$ \cite[1.2.69]{GKOS}). We call \defstyle{monad} of $\tilde x$ the set $\monad(\tilde x)=\{\tilde y\in\GenR^d: \tilde y \approx \tilde x\}$. We denote $\ns(\Omega)=\bigcup_{x\in\Omega} \monad(x)$.
\end{df}
\begin{cor}\label{cor_Geninfty_charac}
Let $u\in\Gen(\Omega)$, $x_0\in\Omega$ and $\tilde x_0\in\widetilde\Omega_c$.
\begin{enumerate}
\item $u\in\widetilde\Gen^\infty_{\widetilde x_0}$ iff there exists $V\in\widetilde{\mathcal V}_{\widetilde x_0}$ such that $u\in\Gentdinfty(V)$
\item $u\in\dot\Gen^\infty_{x_0}$ iff $u\in\Gentdinfty(\monad(x_0))$
\item $u\in\Gen^\infty_{x_0}$ iff there exists $V\in\mathcal V_{x_0}$ such that $u\in\Gentdinfty(\widetilde V)$
\item $u\in\dot\Gen^\infty(\Omega)$ iff $u\in\Gentdinfty(\ns(\Omega))$
\item $u\in\Gen^\infty(\Omega)$ iff $u\in\Gentdinfty(\widetilde\Omega_c)$.
\end{enumerate}
\end{cor}
\begin{proof}
(1) Since $V$ contains an internal sharp neighbourhood $A=\{\tilde x\in\GenR^d: \abs{\tilde x-\tilde x_0}\le\caninf^m\}$ (for some $m\in\N$), this follows from proposition \ref{prop_geninftyinternal} with $A_n=A$, $\forall n$.\\
(2) Since $\monad(x_0)=\bigcap_{n\in\N} A_n$, with $A_n=\{\tilde x\in\GenR^d: \abs{\tilde x - x_0} \le 1/n\}$ internal, this follows from proposition \ref{prop_geninftyinternal}.\\
(3) (see also \cite{OPS}) Since $\widetilde V=[(V)_\eps]$ is internal, this follows from proposition \ref{prop_geninftyinternal} with $A_n=\widetilde V$, $\forall n$.\\
(4) Immediate by (2).\\
(5) (see also \cite{OPS}) Since $\widetilde\Omega=[(\Omega)_\eps]$ is internal, this follows from proposition \ref{prop_geninftyinternal}.
\end{proof}

We now proceed to show that $\Gen^\infty(\Omega)=\dot\Gen^\infty(\Omega)$. For this purpose, we turn to a more quantitative version of an argument used to characterize $\Null(\Omega)$ \cite[Thm.~1.2.3]{GKOS}. 
\begin{prop}\label{prop_Landau}
Let $u=[(u_\eps)_\eps]\in\Gen(\Omega)$ and let $x_0\in \Omega$. Let for each $k\in\N$, $a_k=\sup\{a\in\R: (\forall \beta\in\N^d$ with $\abs\beta=k)(\forall (x_\eps)_\eps\to x_0) (\abs[]{\partial^\beta u_\eps(x_\eps)}\le \eps^a$, for small $\eps)\}$. Let $m\in\N$. If $a_{m+1} < a_m$, then $-a_k$ is a convex function of $k\in\{n\in\N: n\ge m\}$.
\end{prop}
\begin{proof}
(i) Let first $a_1<a_0$. Suppose that $a_1 < \frac{a_0 + a_2}{2}$.
Let $\delta = \min\{\frac{a_0 - a_1}{2}, \frac{a_0 + a_2}{2} -  a_1\}>0$. Let $(x_\eps)_\eps\to x_0$ arbitrary.
Let $(e_i)_{i=1,\dots,d}$ denote the standard basis of $\R^d$. Let $y_\eps = x_\eps + \eps^{a_0-a_1-\delta}e_i$. Then $(y_\eps)_\eps\to x_0$ and $y_\eps\in\Omega$ for small $\eps$. By Taylor's formula,
\[\partial_i u_\eps(x_\eps)=\eps^{a_1-a_0+\delta} (u_\eps(y_\eps)-u_\eps(x_\eps)) - \frac{\eps^{a_0-a_1-\delta}}{2}\partial_i^2 u_\eps(x_\eps + \theta_\eps (y_\eps-x_\eps)),\]
for some $0\le \theta_\eps\le 1$. Since also $(x_\eps + \theta_\eps(y_\eps-x_\eps))_\eps\to x_0$,
\[
\abs{\partial_i u_\eps(x_\eps)}\le \eps^{a_1-a_0 + \delta} \eps^{a_0-\delta/2} + \eps^{a_0-a_1-\delta} \eps^{a_2-\delta/2} \le 2\eps^{a_1+\delta/2}
\]
for small $\eps$. This contradicts the definition of $a_1$. Thus $a_1\ge \frac{a_0 + a_2}{2}$. In particular, $a_2< a_1$.\\
(ii) If $m\in\N$ and $a_{m+1}<a_m$, the same reasoning applied to all $\partial^\alpha u$ with $\abs\alpha=m$ instead of $u$, yields $a_{m+1}\ge\frac{a_m + a_{m+2}}{2}$ (and in particular $a_{m+2}<a_{m+1}$).
\end{proof}
\begin{cor}\label{cor_Landau}
Let $u=[(u_\eps)_\eps]\in\Gen(\Omega)$ and $x_0\in\Omega$. Let $a_k$ ($k\in\N$) as in the previous proposition. Then $u\in\dot\Gen^\infty_{x_0}$ iff $a_k$ is a non-decreasing function of $k\in\N$ iff $\sup\{\sharpnorm{\partial^\beta u(\tilde x)}: \abs\beta=k$ and $\tilde x \approx x_0\}$ is a non-increasing function of $k\in\N$.
\end{cor}
\begin{proof}
If $a_k$ is a non-decreasing function of $k$, then clearly $u\in\Gentdinfty(\monad(x_0))=\dot\Gen^\infty_{x_0}$. Conversely, if there exist $k<l$ such that $a_l<a_k$, then also $a_{m+1}<a_m$ for some $m\in\N$ (with $k \le m < l$). By proposition \ref{prop_Landau}, $a_{m+j}\le a_m - j(a_m - a_{m+1})$, for each $j\in\N$. In particular, $a_j\to-\infty$ as $j\to\infty$. Hence $u\notin\dot\Gen^\infty_{x_0}$.\\
The second equivalence follows from the fact that $a_k=\sup\{a\in\R:$ $(\forall\beta\in\N^d$ with $\abs\beta=k)$ $(\forall \tilde x\approx x_0) (\val(\partial^\beta u_\eps(x_\eps))\ge a)\}$ $=\inf\{\val(\partial^\beta u(\tilde x)): \abs\beta=k$ and $\tilde x \approx x_0\}$, $\forall k\in\N$ (here $\val$ denotes the valuation on $\GenC$).
\end{proof}

\begin{thm}\label{thm_Gendotinfty_equals_Geninfty_globally}
$\Gen^\infty(\Omega)=\dot\Gen^\infty(\Omega)$.
\end{thm}
\begin{proof}
Let $u\in\Gen(\Omega)$.  If $u\notin\Gen^\infty(\Omega)$, we find $K\csub\Omega$,
$\alpha_n\in\N^d$ (for each $n\in\N$), $\eps_{n,m}\in(0,1/m)$ (for each $n,m\in\N$) and $x_{\eps_{n,m}}\in K$ such that $\abs[]{\partial^{\alpha_n}u_{\eps_{n,m}}(x_{\eps_{n,m}})} > \eps_{n,m}^{-n}$, $\forall n,m\in\N$.
Let $L\csub\Omega$ with $K$ contained in the interior of $L$. As $u\in\Gen(\Omega)$, there exists $N\in\N$ such that $\sup_{x\in L}\abs{u_\eps(x)}\le\eps^{-N}$, for small $\eps$. For a fixed $n>N$, we find by the compactness of $K$ some $x_0\in K$ and a subsequence of $(x_{\eps_{n,m}})_{m\in\N}$ converging to $x_0$. With the notations of proposition \ref{prop_Landau}, this implies that $a_0\ge -N$ and $a_{\abs{\alpha_n}}\le -n < -N$. By corollary \ref{cor_Landau}, $u\notin\dot\Gen^\infty_{x_0}$. Hence $u\notin\dot\Gen^\infty(\Omega)$.
The converse inclusion is immediate.
\end{proof}
\begin{cor}
Let $u\in\Gen(\Omega)$. Then $u\in\Gen^\infty(\Omega)$ iff for each $x\in\Omega$, $\sup\{\sharpnorm{\partial^\beta u(\tilde y)}: \abs\beta=k$ and $\tilde y \approx x\}$ is a non-increasing function of $k\in\N$.
\end{cor}
Theorem \ref{thm_Gendotinfty_equals_Geninfty_globally} apparently contradicts a counterexample in \cite{DPS,OPS}. Indeed, the function constructed there in fact belongs to $\Gen^\infty(\R)$, contrary to what is claimed there. Nevertheless, the argument in \cite{DPS,OPS} shows that $\Gen^\infty_{x_0}\subsetneqq \dot\Gen^\infty_{x_0}$ (for $x_0\in\Omega$) and $\widetilde\Gen^\infty_{\widetilde x_0}\subsetneqq \Gentdinfty_{\widetilde x_0}$ (for $\widetilde x_0\in\widetilde\Omega_c$). The latter means that, for $u\in\Gen(\Omega)$, the set of $\tilde x\in\widetilde\Omega_c$ for which $u\in\Gentdinfty_{\widetilde x}$ is not necessarily open in the sharp topology. We now give a corrected version of the counterexample in \cite{DPS,OPS}.
\begin{ex}
Let $\phi\in\Cnt[\infty](\R)$ with $\supp\phi\subseteq[-1,1]$ and $D^k\phi(0)\ne 0$, for infinitely many $k\in\N$. Let $(a_n)_{n\in\N}$ be a strictly decreasing sequence of real numbers tending to $0$. For every $m\in\N$, $m\ge 1$, let $(\eps_{m,n})_{n\in\N}$ be a strictly decreasing sequence of numbers in $(0,1)$ such that $\lim_{n\to\infty}\eps_{m,n}=0$ and
\[\{\eps_{m,n}: n\in\N\}\cap\{\eps_{m',n}: n\in\N\}=\emptyset, \quad\forall m'\ne m.\]
Let $u_\eps=0$ if $\eps\notin\{\eps_{m,n}: n,m\in\N, m\ge 1\}$ and
\begin{align*}
u_{\eps_{m,n}}(x)&=\frac{1}{\eps_{m,n}^{m+1}} \int_0^x \frac{(x-t)^{m-1}}{(m-1)!} \phi\Big(\frac{t-a_m}{\eps_{m,n}^{m+1}}\Big)\,dt, \qquad\forall n,m\in\N, m\ge 1, x\in\R.
\end{align*}
Let $\tilde x_0=[(x_\eps)_\eps]\in\GenR_c$, where
\[
\begin{cases}
x_{\eps_{m,n}}= a_m + \eps_{m,n}^{m+1},& \forall m,n\in\N, m\ge 1\\
x_\eps=0& \forall \eps\notin\{\eps_{m,n}: n,m\in\N, n\ge 1\}.
\end{cases}
\]
Then $u=[(u_\eps)_\eps]\in\Gen(\R)$, $u\in\dot\Gen^\infty_0\setminus\Gen^\infty_0$, $u\in\Gentdinfty_{\widetilde x_0}\setminus\widetilde\Gen^\infty_{\widetilde x_0}$ and $u\notin\Gentdinfty_{a_m}$, for each $m\in\N$.
\end{ex}
\begin{proof}
Let $R\in\R^+$. Since
\[
\sup_{\abs x\le R} \abs{D^k u_{\eps_{m,n}}(x)}
\le
\begin{cases}
\displaystyle\frac{R^{m-k-1}}{(m-k-1)!}\int_\R \eps_{m,n}^{-m-1}\abs[\Big]{\phi\Big(\frac{t-a_m}{\eps_{m,n}^{m+1}}\Big)}\,dt\le e^R\int_\R\abs\phi, & k<m\\
\eps_{m,n}^{-(m+1)(k-m+1)}\sup_{x\in\R}\abs{D^{k-m}\phi(x)}, & k\ge m,
\end{cases}
\]
$u\in\Gen(\R)$.\\
If $k\ge m$, $D^k u_{\eps_{m,n}}(x)=0$ for $x\le a_m-\eps_{m,n}^{m+1}$. Hence for each $k\in\N$, $\abs{D^k u_\eps(x)}\le e^{a_k}\int_\R\abs\phi$ for $\abs{x}\le a_k/2$ and $\eps$ sufficiently small. Thus $u\in\dot\Gen^\infty_0$.\\
Similarly, if $k\ge m$, $D^k u_{\eps_{m,n}}(x_{\eps_{m,n}})=0$, $\forall n\in\N$. Hence $u\in\Gentdinfty_{\widetilde x_0}$.\\
If $k\ge m$, $D^k u_{\eps_{m,n}}(a_m) = \eps_{m,n}^{-(m+1)(k-m+1)} D^{k-m}\phi(0)$, $\forall n\in\N$. Hence $u\notin\Gentdinfty_{a_m}$, for each $m\in\N$. As $a_m\to 0$, this implies also $u\notin\Gen^\infty_0$.\\
Fix $m_0\in\N$ and let $y_{\eps_{m_0,n}}= a_{m_0}$, $\forall n\in\N$ and $y_{\eps}=x_\eps$, if $\eps\notin\{\eps_{m_0,n}: n\in\N\}$. Similarly, for $\tilde y=[(y_\eps)_\eps]$, $u\notin\Gentdinfty_{\widetilde y}$. Yet $\abs{\tilde y-\tilde x_0}\le\caninf^{m_0+1}$. As $m_0\in\N$ arbitrary, $u\notin\widetilde\Gen^\infty_{\widetilde x_0}$.
\end{proof}

\begin{prop}\label{prop_Taylor}
Let $u\in\Gen(\Omega)$ and $\tilde x_0\in\widetilde\Omega_c$.
\begin{enumerate}
\item $u\in\Gentdinfty_{\widetilde x_0}$ iff $(\exists v\in\Gen_c^\infty(\R^d))$ $(\forall \alpha\in\N^d)$ $(\partial^\alpha u(\tilde x_0) = \partial^\alpha v(\tilde x_0))$.
\item $u\in\widetilde\Gen^\infty_{\widetilde x_0}$ iff $(\exists v\in\Gen_c^\infty(\R^d))$ $(\exists V\in\widetilde{\mathcal V}_{\widetilde x_0})$ $(\forall \tilde x\in V)$ $(u(\tilde x) = v(\tilde x))$.
\end{enumerate}
\end{prop}
\begin{proof}
(1) Let $u\in\Gentdinfty_{\widetilde x_0}$ and $\tilde x_0=[(x_{0,\eps})_\eps]$. For each $\alpha\in\N^d$, let $(c_{\alpha,\eps})_\eps$ be a representative of $\partial^\alpha u(\tilde x_0)$ with $\abs{c_{\alpha,\eps}}\le\eps^{-N}$, $\forall\eps$ (this is possible for some $N\in\N$ not depending on $\alpha$). We now let $w_\eps(x)=\sum_{\abs\alpha\le m_\eps}\frac{c_{\alpha,\eps}}{\alpha!}(x-x_{0,\eps})^\alpha$, where $\lim_{\eps\to 0}m_\eps=\infty$. For each $\beta\in\N^d$ and $M\in\N$,
\begin{multline*}
\sup_{\abs x\le M}\abs{\partial^\beta w_\eps(x)}\le \eps^{-N}\sup_{\abs x\le M}\sum_{\abs\alpha\le m_\eps}\frac{\abs[]{\partial^\beta (x-x_{0,\eps})^\alpha}}{\alpha!}\\
\le\eps^{-N}\sup_{\abs x\le M}\sum_{\abs\alpha\le m_\eps}\frac{\abs[]{(x-x_{0,\eps})^\alpha}}{\alpha!}\le \eps^{-N}\sup_{\abs x\le M} e^{d\abs{x-{x_0,\eps}}}\le\eps^{-N-1},
\end{multline*}
for small $\eps$. Hence $w\in\Gen^\infty(\R^d)$. Further, $\partial^\beta w_\eps(x_{0,\eps})=c_{\beta,\eps}$, for small $\eps$, hence $\partial^\beta w(\tilde x_0)=\partial^\beta u(\tilde x_0)$. Choosing $\phi\in\test(B(0,1))$ with $\phi(x)=1$ for $\abs x\le 1/2$, $v(x):=w(x)\phi(x-\tilde x_0)\in\Gen^\infty_c(\R^d)$ and $\partial^\beta v (\tilde x_0)=\partial^\beta w(\tilde x_0)$, $\forall\beta\in\N^d$.\\
(2) If moreover $u\in\widetilde\Gen^\infty_{\widetilde x_0}$, then $f:=u-v\in\widetilde\Gen^\infty_{\widetilde x_0}$ and $\partial^\alpha f(\tilde x_0)=0$, $\forall\alpha\in\N^d$. By proposition \ref{prop_geninftyinternal}, there exist $n,N\in\N\setminus\{0\}$ such that for each $m\in\N$, by the Taylor expansion up to order $m$,
\begin{multline*}
\sup_{\abs {x-x_{0,\eps}}\le\eps^n}\abs{f_\eps(x)}\le\nu_\eps + \sup_{\abs {x-x_{0,\eps}}\le\eps^n}\abs {x-x_{0,\eps}}^{m+1}\sum_{\abs \alpha=m+1}\sup_{\abs {x-x_{0,\eps}}\le\eps^n}\abs{\partial^\alpha f_\eps(x)}\\\le\eps^{n(m+1)-N-1},
\end{multline*}
for small $\eps$ and for some $(\nu_\eps)_\eps\in\Null_\R$. Since $m\in\N$ arbitrary, $f(\tilde x)=0$, i.e., $u(\tilde x)=v(\tilde x)$ for each $\tilde x\in\widetilde\Omega_c$ with $\abs{\tilde x- \tilde x_0}\le\caninf^n$.\\
For the converse implication, we may suppose that $V$ is open in the sharp topology. Let $\tilde x=[(x_\eps)_\eps]\in V$ and $n\in\N$ sufficiently large. Then $u(\tilde y)=v(\tilde y)$, for each $\tilde y\in \widetilde\Omega_c$ with $\abs{\tilde y- \tilde x}\le\caninf^n$. By contraposition, it follows that $(\sup_{\abs{y-x_{\eps}}\le\eps^n}\abs{u_\eps(y)-v_\eps(y)})_\eps\in\Null_\R$. The same argument as in \cite[Prop.~1.2.3]{GKOS} then yields $\partial^\alpha u(\tilde x) = \partial^\alpha v(\tilde x)$, $\forall\alpha\in\N^d$.
\end{proof}

Let $u\in\Gen(\Omega)$, $\tilde x\in\widetilde\Omega_c$ and $\tilde\xi\in {(\R^d\setminus\{0\})}\sptilde_c$. We consider the following microlocal regularity condition on $u$ (on the scale of the sharp neighbourhoods):\\
there exists $v\in\Gen_c(\R^d)$ and a sharp conical neighbourhood $\Gamma$ of $\tilde\xi$ such that
\begin{equation}\label{eqn_microlocal_reg}
\left\{\begin{array}{l}
(\exists V \in\widetilde{\mathcal V}_{\widetilde x}) (\forall \tilde y\in V) (u(\tilde y)= v(\tilde y))\\
\widehat v(\tilde \eta)=0, \quad\forall \tilde \eta\in\GenR^d_{fs}\cap\Gamma.
\end{array}\right.
\end{equation}

\begin{cor}
Let $u\in\Gen(\Omega)$ and $\tilde x\in\widetilde\Omega_c$. If $u\in\widetilde\Gen^\infty_{\widetilde x}$, then $u$ satisfies condition (\ref{eqn_microlocal_reg}) for each $\tilde\xi\in(\GenR^d\setminus\{0\})\sptilde_c$.
\end{cor}
\begin{proof}
If $u\in\widetilde\Gen^\infty_{\widetilde x}$, by proposition \ref{prop_Taylor}, we find $v\in\Gen_c^\infty(\R^d)\subseteq\Gen_\Schwartz^\infty(\R^d)$ coinciding with $u$ on a sharp neighbourhood. By theorem \ref{thm_slowscale_spectrum}, $\widehat v(\tilde\xi)=0$, for each $\tilde\xi\in\GenR^d_{fs}$.
\end{proof}

We do not know if the converse of this corollary is true.
\begin{prop}
If $u\in\Gen(\Omega)$ is $\Gen^\infty$-microlocally regular (cf.\ \cite{GH05}) at $(x_0,\xi_0)\in \Omega\times (\R^d\setminus\{0\})$, then there exist $V\in{\mathcal V}_{x_0}$ and $W\in{\mathcal V}_{\xi_0}$ such that $u$ satisfies condition (\ref{eqn_microlocal_reg}) for each $\tilde x\in\widetilde V=[(V)_\eps]$ and for each $\tilde\xi\in\widetilde W=[(W)_\eps]$.
\end{prop}
\begin{proof}
W.l.o.g., $\abs{\xi_0}=1$. By \cite[Prop.~6.1.3]{GarettoPhD}, there exists $U\in\mathcal V_{x_0}$ and a conic neighbourhood $\Gamma$ of $\xi_0$ such that for each $\phi\in\test(U)$, in particular for some $\phi\in\test(U)$ with $\phi=1$ on a (nongeneralized) neighbourhood of $x_0$,
\[(\exists N\in\N) (\forall m\in\N) \big(\sup_{\xi\in\Gamma}\left<\xi\right>^m\abs[\big]{\widehat{\phi u_\eps}(\xi)}\le\eps^{-N}, \text{ for small }\eps\big).\]
Hence there exists $V\in\mathcal V_{x_0}$ such that for $v=\phi u$ and $\tilde x \in \widetilde V$, $v(\tilde x)=u(\tilde x)$. Further, for $\tilde t\in\GenR_{fs}$ and $\tilde\xi\in\widetilde W$, for some $W\in {\mathcal V}_{\xi_0}$, $\tilde t\tilde\xi\in\widetilde\Gamma$. Hence
$\abs[\big]{\widehat v(\tilde t\tilde\xi)}\le\caninf^{-N}\langle \tilde t\tilde\xi\rangle^{-m}$, for each $m\in\N$. Since $\abs[]{\tilde t\tilde\xi}\ge \caninf^{-a}$, for some $a\in\R^+$, $\langle \tilde t\tilde\xi\rangle^{-m}\le\caninf^{am}$, and $\abs[\big]{\widehat v(\tilde t\tilde\xi)}\le\caninf^{am-N}$, for each $m\in\N$, i.e., $\widehat v(\tilde t\tilde\xi)=0$.
\end{proof}

\end{document}